\documentclass[12pt]{article}
\usepackage[centertags]{amsmath}     
\usepackage{amssymb}
\usepackage{amsthm}
\usepackage{verbatim}                
\usepackage[dvips]{graphics}
\usepackage{enumerate}

\newtheorem{proposition}{Proposition}[section]
\newtheorem{lemma}[proposition]{Lemma}

\newtheorem{conjecture}[proposition]{Conjecture}
\newtheorem{corollary}[proposition]{Corollary}
\newtheorem{theorem}[proposition]{Theorem}
\newtheorem*{thm1.2'}{Theorem 1.2$'$}

\newcommand{\remark}{\noindent\textbf{Remark}}

\newcommand{\dif}{\mathrm{d}}

\newcommand{\vphi}{\varphi}

\newcommand{\NN}{\mathbb{N}}
\newcommand{\RR}{\mathbb{R}}
\newcommand{\TT}{\mathbb{T}}
\newcommand{\ZZ}{\mathbb{Z}}

\newcommand{\Fl}{\mathcal{F}}
\newcommand{\Ll}{\mathcal{L}}

\newcommand{\Vl}{\mathcal{V}}

\newcommand{\DIV}{\mathrm{div}}

\newcommand{\Area}{\mathrm{Area}}
\newcommand{\Vol}{\mathrm{Vol}}
\newcommand{\vol}{\mathrm{vol}}

\newcommand{\spn}{\mathrm{span}}  

\newcommand{\genus}{\mathrm{genus}}

\newcommand{\dd}[1]{{\textstyle\frac{\partial}{\partial #1}}}
\newcommand{\mdd}[1] {{\frac{\partial}{\partial #1}}}
\newcommand{\tdd}[1] {{\textstyle\partial/\partial #1}}

\begin{document}

\title{Spectral Flexibility of \\
Symplectic Manifolds $T^2\times M$}
\author{Dan Mangoubi \footnote{E-mail address: mangoubi@techunix.technion.ac.il} \\
\textit{\small Department of Mathematics, The Technion, Haifa 32000, Israel}}
\date{}
\maketitle

\abstract
 We consider Riemannian metrics compatible with
 the natural symplectic structure on $T^2\times M$,
 where $T^2$ is a symplectic $2$-Torus and $M$ is a
closed symplectic manifold. To each such metric we attach the
corresponding Laplacian and consider its first positive eigenvalue
$\lambda_1$. We show that $\lambda_1$ can be made arbitrarily
large by deforming the metric structure, keeping the symplectic
structure fixed. The conjecture is that the same is true for any
symplectic manifold of dimension $\geq 4$. We reduce the general
conjecture to a purely symplectic question.
\vspace{1.5ex}

\vspace{1ex}

\section{Introduction and Main Results}
Our paper concerns a rigidity result for deformations of
quasi-K\"ahler structures. We begin by giving a historical
overview of the problem considered.  The following theorem due 
to Joseph Hersch is a first one in a
series of rigidity results concerning $\lambda_1$.
\begin{theorem}[\cite{hersch}]
\label{thm:hersch}
 Let $(S^2, g)$ be the $2$-sphere equipped with a
Riemannian metric
 $g$. Then,
  $$\lambda_1(S^2, g) \Area(S^2, g) \leq 8\pi ,$$
  where $\lambda_1(S^2, g)$ is 
the first positive eigenvalue of the Laplacian on $(S^2, g)$.
\end{theorem}
\noindent It is also known that equality occurs if and only if
$(S^2, g)$ is the standard round sphere.
Hersch's theorem was extended to any surface:
\begin{theorem}[\cite{yangyau}]
\label{thm:yangyau} Let $(\Sigma, g)$ be a closed Riemannian
surface. Then
  $$\lambda_1(g) \Area(\Sigma, g) \leq 8\pi(\genus(\Sigma)+1).$$
\end{theorem}
\noindent Here, the constant on the right hand side is not
optimal.

 The next result shows non-rigidity in higher dimensions.
\begin{theorem}[\cite{col-dod}]
Let $M$ be a closed Riemannian manifold of dimension $>2$. Then
one can find a Riemannian metric of total volume~$1$ and of
arbitrarily large $\lambda_1$.
\end{theorem}
 If one fixes a conformal class on a manifold $M$, then one
  recovers rigidity for $\lambda_1$:
 \begin{theorem}[\cite{soufi-ilias, fried-nadi}]
 Let $(M, g)$ be a closed Riemannian manifold of dimension $n$. Then
 $$\lambda_1 (fg) \Vol(M, fg)^{2/n} \leq C(g), $$
 where $f$ is any positive function on $M$ and $C(g)$ is a constant
 independent of $f$.
 \end{theorem}
\noindent In particular, by the Uniformization Theorem for the
Riemann sphere this theorem extends Theorem~\ref{thm:hersch}.

In~\cite{pol-sfev} Polterovich considers Theorem~\ref{thm:yangyau}
in the K\"ahler and quasi-K\"ahler categories. In order to
explain this, we recall the
definitions of K\"ahler and quasi-K\"ahler structures: Let $(M,
\omega)$ be a closed symplectic manifold.
 Let $J$ be an almost complex structure on $M$ which is
 compatible with $\omega$, i.e.~$\omega(J v, J w) = \omega(v, w)$
 and $\omega(v, J v)>0$
   unless $v=0$.
   Let $g$ be the corresponding Riemannian metric
   $g(v, w) :=\omega(v, J w)$.
   The quadruple $(M, \omega, J, g)$ is called a
   \emph{quasi-K\"ahler} structure on $M$. If $J$ is a complex
   structure, it is called a \emph{K\"ahler} structure.

Having these terms in mind we reformulate
Theorem~\ref{thm:yangyau} as
\begin{thm1.2'}[\cite{yangyau}]
  Let $(\Sigma, \omega)$ be a closed symplectic surface.
  For any quasi-K\"ahler structure $(\Sigma, \omega, J, g)$,
  one has $$\lambda_1(g)\Area(\Sigma, \omega)\leq 8\pi(\genus(\Sigma)+1).$$
\end{thm1.2'}

The equivalence of the two formulations follows from 
Moser's argument~\cite[Th.~3.17]{mcd-sal}. 
This argument shows that fixing the total area of a closed 
surface is equivalent to fixing the symplectic structure
$\omega$ on it. 
Indeed,
according to this argument, if $(\omega_t)$, $0\leq t\leq 1$, is a smooth family
of symplectic forms on a compact manifold, all in the same cohomolgy class,
then there exists a flow $(\vphi_t)$ for which $\vphi_t^*(\omega_t) = \omega_0$.
Now, for any two symplectic forms $\omega_0$ and $\omega_1$,
we define $\omega_t=(1-t)\omega_0 + t\omega_1$.
In dimension two $(\omega_t)$ is a family of symplectic forms,
since this is a family of area forms. In addition, if $\omega_0$ and $\omega_1$
have the same total area, then all the $\omega_t$'s are cohomologous. 
Thus, the conditions of Moser's theorem are satisfied and we conclude
that $\omega_0$ and $\omega_1$ are symplecto\-morphic.
 
Also, linear algebra in $\RR^2$ shows
that any Riemannian metric on a symplectic surface comes from a
quasi-K\"ahler structure whose quasi-K\"ahler form is $\omega$
(moreover, any almost complex structure on a closed surface is
integrable).

The following theorem is proved in~\cite{pol-sfev}. Its proof is based
on results from \cite{bly} and \cite{dema}.
\begin{theorem}[\cite{pol-sfev}]
Let $(M, \omega)$ be a closed symplectic manifold. Suppose
$\omega$ is a rational form.
Let $g$ be a K\"ahler metric whose K\"ahler form is $\omega$.
Then
$$\lambda_1(g) \leq C(\omega),$$
  where $C(\omega)$ is independent of $g$.
\end{theorem}
\noindent It is still an open question whether the same is true
for any real symplectic form $\omega$.
On the other hand, if we consider quasi-K\"ahler metrics, then the
conjecture is
\begin{conjecture}
\label{conj:general}
Let $(M, \omega)$ be a closed symplectic manifold of dimension
$\geq 4$. Then, there exists a quasi-K\"ahler structure on it with
arbitrarily large $\lambda_1$.
\end{conjecture}

For $\TT^4\times M$ we have
\begin{theorem}[\cite{pol-sfev}]
Let $(\TT^4, \sigma)$ be the standard symplectic $4$-torus.
Let $(M, \omega)$ be a closed symplectic manifold.
Then, on $(\TT^4\times M, \sigma\oplus\omega)$ there exists
a quasi-K\"ahler structure with
arbitrarily large $\lambda_1$.
\end{theorem}

In his proof Polterovich shows existence of certain plane
distributions on $(\TT^4\times M, \sigma\oplus\omega)$ along
which one can deform an almost complex structure $J$.

 The main novelty of our paper is the introduction of non-regular
distributions in this procedure. These distributions do not have a
constant dimension. Their dimension may drop on exceptional subsets.
 The precise definition is given in Section~\ref{sec:overview}.
  The advantage of our
method is in two aspects: First, singular distributions may exist
on manifolds on which regular distributions do not exist simply
due to topological restrictions. Thus, we extend the family of
manifolds for which Conjecture~\ref{conj:general} is true. In
particular, we prove
 \begin{theorem}
 \label{t2xm}
 Let $(\TT^2, \sigma)$ be the standard symplectic $2$-torus.
Let $(M, \omega)$ be a closed symplectic manifold.
Then, on $(\TT^2\times M, \sigma\oplus\omega)$ there exists
a quasi-K\"ahler structure with
arbitrarily large $\lambda_1$.
\end{theorem}
 Second, our method shows
 that Conjecture~\ref{conj:general} follows from
a purely symplectic conjecture. Namely,
\begin{conjecture}
\label{conj:total-symp}
Let $(M, \omega)$ be a closed symplectic manifold of dimension
$n\geq 4$. Then one can find on $(M, \omega)$  an isotropic singular distribution
$\Ll$ which satisfies H\"ormander's condition.
\end{conjecture}
\noindent Singular distributions and H\"ormander's condition are
discussed in Section~\ref{sec:overview}. We prove
\begin{theorem}
\label{thm:symp-spectral}
  Conjecture~\ref{conj:general} follows from Conjecture~\ref{conj:total-symp}.
\end{theorem}
We would like to emphasize that any answer to
Conjecture~\ref{conj:total-symp} would be very interesting. As
mentioned, a \emph{positive} answer will resolve
Conjecture~\ref{conj:general}, while a \emph{negative} answer will
presumably lead to a new type of symplectic rigidity.

\subsection{Acknowledgements}
I would like to express my deep gratitude and
appreciation to my advisor
Leonid Polterovich for his support and encouragement.

The introduction of this paper was improved according to the remarks
of an anonymous referee. I am grateful to him.

This paper was written while being a student at the Technion, Haifa, Israel.
I would like to thank the Technion for its generous support.
\section{Overview of the Main Ideas}
 \label{sec:overview}
The basic object we use is a finite set $X=\{X_k\}_{k=1}^{N}$ of
smooth vector fields, which satisfy
\begin{itemize}
\item[1.] $\exists A\in\NN$ such that the dimension of
$$\Ll|_{x}=\spn(X_1(x), \ldots, X_N(x))$$
is $A$ almost everywhere. The complement of $\{\dim \Ll|_{x}=A\}$
is called the singular set.
\item[2.] The iterated commutators of the $X_k$'s and which are of lengths
$\leq r$ span $T_x M$ for every $x$.
\end{itemize}
A set $X$ which satisfies property~1 will be called a
\emph{singular distribution}, and will be denoted by $\Ll$.
Property~2 is called the \emph{$r$-H\"or\-mander condition.}

 The idea to use H\"or\-mander regular distributions
 in the context of $\lambda_1$ can be found in \cite{bour-ber}
 and \cite{pol-sfev}. The new idea we introduce is to allow also
 distributions which might not have constant dimension.
 The dimension of a singular distribution might drop on subsets of measure~$0$.
 Singular distributions are subject to less topological restrictions, and
thus can be used where smooth distributions do not exist.

To a H\"or\-mander distribution associated with~$X$, we can attach
the hypo\-elliptic operator
$$L= \sum X_j^{*} X_j\, ,$$
where $X_j^{*} = -X_j-\DIV(X_j)$. We start with an arbitrary
compatible Riemannian metric $g$ and we deform it in a way which
is related to the distribution associated with~$X$. We would like
to relate the spectral properties of $L$ to the spectrum of the
Laplacian of a deformed metric. The estimate of the eigenvalues of
the deformed metric is done through the variational principle.

 Unlike in~\cite{pol-sfev}, we can only deform the
 Riemannian metric structure on $M$
 away from the singular set of $\Ll$.
 The new difficulty we face in our estimates is the fact that the test
 functions near the singular set do not feel the deformation,
 and the estimate of~\cite{pol-sfev} does not apply anymore.

 To overcome this problem we apply the theory of an\-iso\-tropic So\-bo\-lev
 spaces as developed in~\cite{roth-stein}, and the machinery
 of fractional Sobolev spaces also known as Bessel Potential Spaces.

Theorem~\ref{t2xm} is a direct consequence of the following two
theorems.

\begin{theorem}
  \label{deform}
  Suppose $\Ll$ is an isotropic singular distribution on $(M, \omega)$
  which satisfies H\"{o}r\-mander's condition.
    Then, one can find a compatible almost complex structure $J$ on $M$
    with arbitrarily large $\lambda_1(M, \omega, J)$.
\end{theorem}

\begin{theorem}
  \label{thm:t2xm}
  Let $(M, \sigma_M)$ be a closed symplectic manifold.
  On $(T^2\times M, \sigma_T\oplus\sigma_M)$ there exists an isotropic
  singular distribution which satisfies H\"{o}r\-mander's condition.
\end{theorem}
\noindent The proofs are given in section~\ref{prf-dfrm} and
section~\ref{prf-t2xm} below.
\section{Proof of Theorem~\ref{deform}}
\label{prf-dfrm}

We assume that $(M, \omega)$ is a closed symplectic manifold of
dim $n$, and $\Ll$ is a singular isotropic distribution which
satisfies the $r$-H\"or\-mander condition. We can find on $M$ an
arbitrary quasi-K\"ahler structure $(M, \omega, J_1,g_1)$.
\subsection{A Deformation of the Almost Complex Structure}
  We introduce a deformation of the almost
  complex structure $J_1$, which keeps it compatible with $\omega$.
  The deformation will be defined with the help of the distribution $\Ll$.

  Let $B$ be the singular set of $\Ll$.
   Let $\phi$ be a positive function on $M$ which is identically $1$
    in a neighborhood of $B$.

    Consider the $g_1$-orthogonal decomposition at a point
    $x\in M$.
    $$T_x M = \Ll|_x\oplus J_1 \Ll|_x \oplus \Vl_x\, . $$
    Set
    \begin{equation}
      J_\phi = \left\{\begin{array}{lll}
                     (1/\phi) J_1 &  & \mbox{on } \Ll \, ,\\
             \phi J_1 &  & \mbox {on } J_1 \Ll \, , \\
             J_1 & & \mbox{on } \Vl \, .
                     \end{array}\right.
    \end{equation}
    One can check that $J_\phi$ defines a quasi-K\"ahler structure
    $(M, \omega, J_\phi, g_\phi)$.

%
\subsection{Estimation of $\lambda_1$}
    We estimate $\lambda_1(M, \omega, J_\phi, g_\phi)$ by using the
    variational formula.

    Let $K$ be such that $g_1(X_j, X_j)\leq K$ at every point.
    Since $\Ll$ satisfies H\"or\-mander's condition, we can
    apply the following version of the Poincar\'{e}--Sobolev
    inequality. Its proof is given in section~\ref{sec:poincare}.
    \begin{theorem}
      If $2-\frac{4}{rn+2}<p<2$,
      then
      $$ \|u-\bar{u}\|_{L^2} \leq
      C_p \sum_{j=1}^{N}\left( \int_M| X_j u|^p\,\dif(\vol)
                                        \right)^{1/p}, $$
  for all $u\in C^\infty(M)$. Here $\bar{u}$ is the average of $u$ on
  $M$, and $C_p$ is independent of $u$.
   \end{theorem}

  Fix $p$ as in the theorem. Let
  $\langle\cdot,\cdot\rangle_{\phi}\,$,
  $\|\cdot\|_{\phi}\,$, $\nabla_{\phi}$ denote the inner product,
  the norm and the gradient respectively of the Riemannian metric
  $g_\phi$. Note that the volume form $\vol_\phi=\vol$ for any
  $\phi$, since $g_\phi$ is compatible with $\omega$.
  For any function $u$ of zero average, we obtain
    \begin{eqnarray*}
      \|u\|_{L^2}\! &\leq&\!
         C_p\sum_{j=1}^{N} \left(\int_M|
         \langle\nabla_{\phi} u, X_j\rangle_{\phi} |^p\,\dif(\vol)
                                   \right)^{1/p} \\
       &\leq&\! C_p N\sqrt{K}
           \left(\int_M \| \nabla_\phi u \|_{\phi}^{p}
           \frac{1}{\phi^{p/2}}\,\dif(\vol)
                                         \right)^{1/p} \\
       &\leq&\! C_p N\sqrt{K}
           \left(\int_M  \| \nabla_\phi u \|_{\phi}^2\,\dif(\vol)
           \right)^{\frac{1}{2}}
         \left(\int_M\frac{1}{\phi^{p/(2-p)}}\,\dif(\vol)\right)^{\frac{2-p}{2p}},
    \end{eqnarray*}
    where we used the fact that $\|X_j\|_{\phi}=
    \phi^{-1/2}\|X_j\|_1\leq \phi^{-1/2}K^{1/2}$.
    Hence,
      $$ \frac{\int_M  \| \nabla_\phi u \|_{\phi}^2\,\dif(\vol)}
      {\int_M |u|^2\,\dif(\vol)}\geq (C_p N)^{-2}K^{-1}
      \|1/\phi\|_{L^{p/(2-p)}}^{-1}\, .$$
    We can choose
    $\phi$ such that
    $\|1/\phi\|_{L^{p/(2-p)}}$
    is very small. By the variational principle we conclude that $\lambda_1$
    can be made arbitrarily large.

\section{Proof of Theorem~\ref{thm:t2xm}}
\label{prf-t2xm}
Let us denote the coordinates on the torus $T^2=\RR^2/(2\pi\ZZ)^2$
by $x, y$ modulo $2\pi$. Let $z$ denote a point on $M$.

\subsection{Construction of a Distribution}
Let $h(z)$ be a Morse function on $M$. Let $Z_1,\ldots, Z_{2d}$ be
vector fields on $M$, which span $T_z M$ for every $z\in M$.
Let
$$\phi_{2j-1}(x)=\sin(j x),\ \phi_{2j}(x) = \cos(j x),\ \ (1\leq j\leq d).$$
\remark \textbf{ and Notation}. If we let $Q(t)= \prod_{j=1}^{d} (t^2+j^2)$,
then $(\phi_l)_{l=1}^{2d}$ is a basis of solutions to the
differential equation $Q(\dd{x})(u) = 0\, .$ We will also write
$Q(t)(u)=0$, where $t$ acts by $\tdd{x}$.

Define the following two vector fields on $T^2\times M$:
\begin{eqnarray}
 V(x, y, z) &:=& \mdd{x} +  h(z)\mdd{y}\, , \\
 W(x, y, z) &:=& \sum_{l=1}^{2d} \phi_l(x) Z_l(z)\, .\nonumber
\end{eqnarray}
We claim that these two vector fields span a singular isotropic
H\"or\-mander distribution $\Ll$.
It is clear that $\Ll$ is isotropic, so it remains to check
that the H\"or\-mander condition
and the singularity property are satisfied for $\Ll$.

\subsection{$\Ll$ is H\"{o}r\-mander}
For two vector fields $X, Y$, set $D_X(Y) := [X,Y]$.

We will prove that the following iterated commutators
$$ V, W, D_V^\alpha(W), D_V^\gamma D_W D_V^\beta(W),$$
where $0\leq\alpha<2d$, $0\leq \beta< 4d$ and $0\leq\gamma< 2d$ span
the tangent space at every point.


\begin{description}
\item[\mdseries\scshape{Step~1.}] We compute by induction on $\alpha$
 $$
   D_V^\alpha(W)(x, y, z) = \sum_{l=1}^{2d} \phi_l^{(\alpha)}(x)Z_l(z)
   - \alpha \phi_l^{(\alpha-1)}(x)(Z_l h)(z)\mdd{y}(y)\, .
 $$
\item[\mdseries\scshape{Step 2.}] We prove that the vector fields
  $\{D_V^\alpha(W)\}_{0\leq \alpha<2d}$ span
  $Z_l$ and $(Z_l h) \tdd{y}$ at every point for ${1\leq l \leq 2d}$.
\end{description}

By Step~1, for any polynomial $F(t) = \sum a_k t^k$ we have

$$F(D_V(W))=\sum_k a_k D_V^k(W) = \sum_{l=1}^{2d}(F(t)\phi_l)Z_l
- \sum_{l=1}^{2d} (F'(t)\phi_l)(Z_l h) \mdd{y}\, .$$
In particular, for $F(t) = t^{r-1} Q(t) = \sum_{k} q_k t^{r+k-1}$ we have
\begin{equation*}
  \xi_r:=\sum_{k=0}^{2d} q_k D_V^{r+k-1}(W) = -\sum_{l=1}^{2d} (t^{r-1} Q'(t)\phi_l) (Z_l h) \mdd{y} \, .
\end{equation*}

\begin{lemma}
  \label{basis}
  $Q'(t)\phi_l = b_l \phi_l'$ for some integer $b_l\neq 0$, where $1\leq l\leq 2d$.
\end{lemma}
We give the proof below. Thus,
\begin{equation*}
  \xi_r = -\sum_{l=1}^{2d}b_l(\phi^{(r)}_{l})(Z_l h) \mdd{y}\, .
\end{equation*}
 Since the $(2d\times 2d)$
Wronskian matrix $(\phi^{(r)}_l)_{rl}$ is non-degenerate and $b_l\neq0$,
it follows that the $\xi_r$'s
$(1\leq r\leq 2d)$ span the vector fields $(Z_l h) \tdd{y}$ at every point.
Finally, from Step~1 we see
that the $Z_l$'s are spanned by vectors of the forms
$D_V^{\alpha}(W)$ and $(Z_l h)\tdd{y}$.
\begin{description}
\item[\mdseries\scshape{Step 3.}] We show that
$$\{D_V^\gamma D_W D_V^\beta(W), D_V^\alpha(W)\}
\mbox{ span } \sum_{k, l=1}^{2d} \phi_k^{(i)}\phi_l^{(j)}
(Z_k Z_l h)\mdd{y} $$
for $0\leq i, j< 2d$, where $0\leq\alpha<2d$, $0\leq\beta<4d$ and $0\leq\gamma< 2d$.
\end{description}
%
%

  We compute:
    $$D_W D_V^\beta(W) = \sum_l a_l^{\beta}Z_l -
  \beta\sum_{k, l} \phi_k\phi_l^{(\beta-1)}(Z_k Z_l h) \mdd{y},$$
  for some $C^\infty$-functions $a_l^{\beta}$ on $T^2\times M$.

  Since we can span the $Z_l$'s at every point of $T^2\times M$
  (by Step~2), it follows that
  we can also span $\sum_{k, l} \phi_k\phi_l^{(j)} (Z_k Z_l h)\tdd{y}$
  for \mbox{$0\leq j \leq 2d$}. Let
  $$Y=\sum_{k, l} \phi_k^{(a)}\phi_l^{(b)} (Z_k Z_l h)\tdd{y}\ .$$
  By repeated uses of the formula
   $$D_V(Y) = \sum_{k, l} \left(\phi_k^{(a+1)}\phi_l^{(b)} +
                  \phi_k^{(a)}\phi_l^{(b+1)}\right)(Z_k Z_l h)\mdd{y}\ ,$$
   we see that the given vectors span
   $\sum_{k, l} \phi_k^{(i)}\phi_l^{(j)} (Z_k Z_l h)\tdd{y}$.
\begin{description}
\item[\mdseries\scshape{Step 4.}] $$\{D_V^\gamma D_W D_V^\beta(W), D_V^\alpha(W)\}
\mbox{ span } (Z_k Z_l h)\mdd{y} $$
for any $1\leq k,l\leq 2d $, where $0\leq\alpha<2d$, $0\leq\beta<4d$ and $0\leq\gamma<2d$.
\end{description}

The $(2d\times 2d)$ matrix $B=(\phi_k^{(i)})_{ik}$
 is non-degenerate at every point.
 Hence the matrix $B\otimes B = (\phi_k^{(i)} \phi_{l}^{(j)})_{ij, kl}$
 is non-degenerate. It follows from Step~3 that
 we can span $(Z_k Z_l h)\tdd{y}$.
\begin{description}
 \item[\mdseries\scshape{Step 5.}] We finish by noting that since $h$ is a Morse function,
 we have $\tdd{y}$ at every point by Step~4 and Step~2.
 We use $V$ to get also $\tdd{x}$.
\end{description}

\subsubsection{Proof of Lemma~\ref{basis}}
  Let $1 \leq j\leq d$.
  $Q(t)=P(t^2)$, where
  $P(t)=\prod_{k=1}^{d}(t+k^2)$. We can write
  $P'(t)=F_j(t)(t+j^2) + b_j$,
  for some polynomial $F_j$ and integer $b_j$. Moreover,
  $b_j\neq 0$, since $P'(-j^2)=\prod_{k\neq j}(-j^2+k^2)\neq 0$.
  Hence, $Q'(t)=2tP'(t^2)=2tF_j(t^2)(t^2+j^2) + 2b_j t$.

  Since   $\phi_{2j-1}$ and $\phi_{2j}$ are solutions for
  the second order differential equation $(t^2+j^2)u = 0$,
  we conclude that $Q'(t)\phi_l = 2b_{\lfloor\frac{l+1}{2}\rfloor} \phi_l'$.

\subsection{$\Ll$ is a Singular Distribution}
We will show that $\dim(\spn(V, W))=2$ almost everywhere.
 It suffices to check that $W$ vanishes on a closed set of measure $0$.

Suppose $W(p)=0$. Choose coordinates $(x, y, z_1, \ldots z_n)$
in a neighborhood of $p$. We can write
  $$Z_l = \sum_{j=1}^{n} a_{jl} \mdd{z_j}$$
for some $C^{\infty}$-functions $a_{jl}$. The $(n \times 2d)$
matrix $A=(a_{jl})$ is of full rank~$n$. When written in
coordinates, the equation $W(p) = 0$ becomes the following system
of equations:
 $$ \sum_{l=1}^{2d} \phi_l(x) a_{jl}(z) = 0,\ (1\leq j\leq n),$$
 which can also be written as
 $ A(z) \phi(x) = 0$, where $\phi$ is the column vector of the $\phi_l$'s.
 Let $E=\{(x,y,z)\in \TT^2\times M \, : A(z) \phi(x) =0\}$.
 It is clear that $E$ is closed.
 For each pair $(y, z)$ let
   $$E_{y,z} = \{ x \, : (x,y,z)\in E\}$$
 be a subset of the circle.

 We show that $E_{y,z}$ is a discrete set, hence finite of measure $0$.
 Let $x\in E_{y,z}$. The matrix $(\phi_l^{(m)}(x))_{lm}$, $(1\leq l \leq 2d, 0\leq m < 2d)$
 is of full rank. Therefore, for some $0<m <2d$, $A(z)\phi^{(m)}(x) \neq 0$.
 From here it follows that $x$ is an isolated point in $E_{y,z}$.

 We have shown that, for every fixed pair $(y,z)$, $E_{y, z}$ is of measure~$0$.
 By Fubini's Theorem it follows that $E$ is of measure~$0$.
\section{Sobolev spaces associated with H\"or\-mander Distributions}
\label{sobolev}
Let $M$ be a closed manifold, equipped with a volume form. In this
section we give an overview of a family of Sobolev spaces
associated with a H\"or\-mander set of vector fields $X$ (see
section~\ref{sec:overview}).
These spaces were defined and studied in~\cite{roth-stein}. 

For $1<p<\infty$, let
$$W^{1, p}_{X}(M) = \{f\in L^p(M)\,:\, \forall j \,X_j f\in L^p(M) \}. $$
The norm on $W^{1,p}_{X}(M)$ is given by
  $$ \|u\|_{1,p,X} = \|u\|_{L^p} + \sum_{j=1}^{N} \|X_j u\|_{L^p}\ .$$
\subsection{Bessel Potential Spaces}
In order to formulate embedding theorems for the spaces $W_X^{1,p}$
we review a family of fractional Sobolev spaces also known as Bessel
potential spaces. For more details on these spaces
see~\cite[ch.~5]{stein}, \cite[pp.~219--222]{adams},
\cite{arons-smith, adams-arons-smith, adams-arons-hanna, calderon, ar-mu-sz, triebelI, triebelII}.

For $s\geq 0$ and $1<p<\infty$, let
 $$W^{s,p}(\RR^n) = \{f\in L^p(\RR^n) \,:\,
 \Fl^{-1}((1+|x|^2)^{s/2} \Fl(f))\in L^p(\RR^n)\},$$
where $\Fl$ is the Fourier transform on tempered distributions.
The norm
 $$ \| f \|_{s,p} := \| \Fl^{-1}((1+|x|^2)^{s/2} \Fl(f)) \|_{L^p} $$
 makes $W^{s,p}(\RR^n)$ into a Banach space. For an integer $s$, $W^{s,p}$ coincides
 with the classical Sobolev space (not true for $p=1$).

\remark. If we set $G_s =\Fl((1+|x|^2)^{-s/2})$ (the Bessel
potential), then $W^{s,p} = G_s * L^p$, and $\|f\|_{s,p}=\|G_{-s}*
f\|_{L^p} $, where $*$ denotes convolution.

Next we define $W^{s,p}(\Omega)$ for a bounded $C^{\infty}$-domain $\Omega \subset \RR^n$
as the space of restrictions from $W^{s,p}(\RR^n)$.
The natural norm here is the quotient norm
$$\|f\|_{s, p, \Omega} = \inf_{\tilde{u}} \|\tilde{u}\|_{s, p, \RR^n},$$
 where $\tilde{u}$ is an extension of $u$ to $\RR^n$.

The space $W^{s,p}(\RR^n)$ is invariant under
$C^\infty$-diffeo\-morphisms of $\RR^n$ and is a
$C^\infty(M)$-module. Therefore we can define $W^{s,p}(M)$ on a
closed manifold $M$ in the standard way (\cite{adams-arons-hanna},
\cite[Ch.~7]{triebelII}). Namely, on $M$ choose a finite partition
of unity $\psi_i$ subordinated to an atlas $(\Omega_i, h_i)$ of
$M$. We define
$$ W^{s,p}(M) = \{f:M\to\RR \,:\, (\psi_i u)\circ h_i \in W^{s,p}(\Omega_i)\}, $$
where a norm is given by
$$\|u\|_{s, p, M} := \sum_i \|(\psi_i u)\circ h_i\|_{s, p, \Omega_i}.$$
\subsection{Embedding Theorems}
Suppose that $X$ satisfies the $r$-H\"{o}r\-mander condition. We
have
\begin{theorem}[{\cite[Theorem~17]{roth-stein}}]
   The space $W^{1, p}_{X}(M)$ is continuously embedded in $W^{1/r,p}(M)$.
\end{theorem}

Let $\dim(M)=n$. The (fractional) Sobolev embedding theorem is
(see e.g.~\cite[pp. 470--471]{arons-smith},
\cite[Theorem~7.63]{adams}, \cite[sec.~2.7]{triebelI}):
\begin{theorem}
 \label{sob}
 If $sp<n$ and $0\leq 1/p-1/q\leq s/n$, then the space
 $W^{s,p}(M)$ is continuously embedded in $L^{q}(M)$.
\end{theorem}
%

The (fractional) Rellich--Kondrachov Theorem for a closed manifold $M$ is
\begin{theorem}
   If $sp<n$ and $1/p-1/q< s/n$, then
  the space $W^{s,p}(M)$ is compactly embedded in $L^q$.
\end{theorem}
\remark. This result is stated in~\cite[sec.~2.5.1]{edm-tri} for
bounded domains in $\RR^n$. The same result for closed manifolds
is deduced in a standard way, and we omit the proof.

\begin{corollary}
 \label{X-sob} 
 If $p/r<n$ and $0\leq 1/p-1/q\leq 1/(rn)$,
 then the space $W^{1,p}_{X}(M)$ is continuously embedded in
 $L^{q}(M)$.
\end{corollary}

\begin{corollary}
 \label{X-rk} 
 If $p/r<n$ and $1/p-1/q<1/(rn)$,
 then the space $W_{X}^{1,p}(M)$ is
 compactly embedded in $L^{q}(M)$.
\end{corollary}
%
\subsection{Poincar\'e--Sobolev Inequalities}
\label{sec:poincare}
 Define the energy corresponding to the anisotropic
Sobolev space $W_X^{1, p}$ by
$$E_{X,p}(u) =\sum_{j=1}^{N} \|X_j u\|_{L^p}\, .$$
\begin{theorem}
  \label{poincare}
  Let $p>1$, $1/p-1/q<1/(rn)$.
  Then
  $$ \|u-\bar{u}\|_{L^q} \leq C_{p,q} E_{X,p}(u),$$
  for any $u\in C^{\infty}(M)$, where $C_{p,q}$ is independent of $u$ and
  where $\bar{u}$ is the average of $u$ on $M$.
\end{theorem}
\begin{proof}
\underline{\scshape{Step 1}}.
We can assume that $q=p$. Indeed, if $q>p$, then by
Corollary~\ref{X-sob} one has
  $$\|u-\bar{u}\|_{L^q} \leq C_{p,q} (\|u\|_{L^p}+E_{X,p}(u))\, .$$
  Joining this with the result for $q=p$ gives the desired
  inequality.

  If $1<q<p$, then we conclude the inequality from H\"older's
  inequality and the result for $q=p$.

\underline{\scshape{Step 2}}.
  It is enough to prove the result for functions with average~$0$.
  Suppose there does not exist a constant $C_{p,p}$ as above.
  Then we can find a sequence $(u_n)$, $u_n \in C^\infty(M)$, with average~$0$,
 such that ${\|u_n\|}_{L^p}=1$ and for all $j$
$X_j u_n \stackrel{L^p}{\to} 0$. Hence, the sequence $(u_n)$ is
bounded in $W^{1,p}_{X}$. By Corollary~\ref{X-rk}, we may assume
that $(u_n)$ converges to $u$ in $L^p$. On the one hand $X_j u_n
\to X_j u$ in distribution sense, and on the other hand $X_j u_n
\stackrel{L^1}{\to} 0$ by H\"{o}lder's inequality, and therefore also
in distribution sense. We obtain that for all $j$ $X_j u = 0$ in
distribution sense.

\underline{\scshape{Step 3}}. Any $Y\in T_x M$ can be expressed as
a sum $Y=\sum a_j Y_j$, where the $Y_j$'s are iterated commutators
of the $X_j$'s. Hence, $Yu=0$ for any smooth vector field $Y$.
Therefore, $u$ is a constant function. Since it has average $0$,
we conclude $u=0$. This is a contradiction to $1=\| u_n \|_{L^p}
\to \| u\|_{L^p}$.
\end{proof}

\providecommand{\bysame}{\leavevmode\hbox to3em{\hrulefill}\thinspace}
\providecommand{\MR}{\relax\ifhmode\unskip\space\fi MR }
\providecommand{\MRhref}[2]{%
  \href{http://www.ams.org/mathscinet-getitem?mr=#1}{#2}
}
\providecommand{\href}[2]{#2}

%

\end{document}